\documentclass[12pt,a4paper]{article}
\usepackage{amsmath}
\usepackage{amsfonts}
\usepackage{amssymb}
\usepackage{makeidx}
\usepackage{fullpage}
\usepackage{graphicx,float}
\usepackage{xcolor,subfigure}
\usepackage[export]{adjustbox}
\definecolor{webgreen}{rgb}{0,.5,0}
\definecolor{webbrown}{rgb}{.6,0,0}
\usepackage[colorlinks=true,
linkcolor=webgreen,
filecolor=webbrown,
citecolor=webgreen]{hyperref}
\newcommand{\seqnum}[1]{\href{http://oeis.org/#1}{\underline{#1}}}
\newcommand{\Figs}[1]{\hyperref[#1]{Figure~\ref*{#1}}}

\title{\bf An approximate Jerusalem square whose side equals a Pell number}
\author{Franck Ramaharo\\
\small D\'epartement de Math\'ematiques et Informatique\\[-0.8ex]
\small Universit\'e d'Antananarivo\\[-0.8ex] 
\small 101 Antananarivo, Madagascar\\
\small\href{mailto:franck.ramaharo@gmail.com}{\tt franck.ramaharo@gmail.com}\\
}

\date{\small\today}
\begin{document}
\maketitle
\begin{abstract}
We take advantage of the properties of the Pell numbers to construct an integer version  of the Jerusalem square fractal.

\end{abstract}

\section{Introduction}
Eric Baird \cite{BairdSquare} first introduced the Jerusalem square in 2011. This fractal object  can be constructed as follows.
\begin{enumerate}
\item Start with a square.
\item Cut a cross through the square so that the corners then consist of four smaller scaled copies, of rank $ +1$,  of the original square, each pair of which being separated by a smaller square, of rank $ +2 $, centered along the edges of the original square. The scaling factor  between the side length of the squares of consecutive rank is constrained to be constant.
\item Repeat the process on the squares of rank $+ 1 $ and $ +2 $, see \Figs{Fig:jerusalem-square}.
\begin{figure}[H]
\centering
\includegraphics[width=\linewidth]{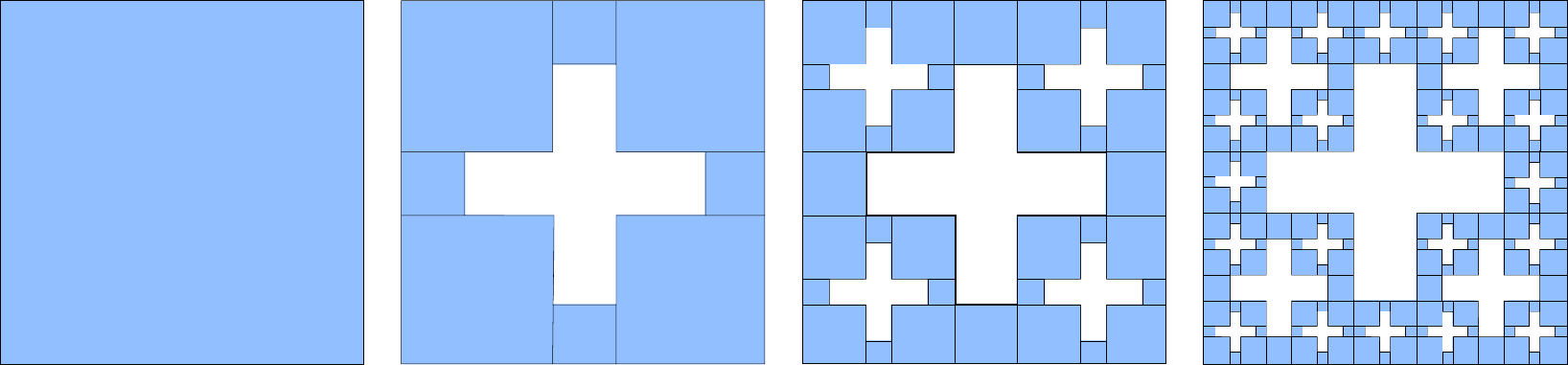}
\caption{The iterative construction of the Jerusalem square.}
\label{Fig:jerusalem-square}
\end{figure}

\end{enumerate}
Let $ n $ be a nonnegative integer, and let $ \ell_n $ denote the side length of the square at the $ n $-th iteration. Then, 
\begin{equation}\label{eq:reccurence}
\ell_n=2\ell_{n+1}+\ell_{n+2}.
\end{equation}
The  scaling factor constraint implies
\begin{equation}\label{eq:ratio}
k=\dfrac{\ell_{n+1}}{\ell_n}=\dfrac{\ell_{n+2}}{\ell_{n+1}}.
\end{equation}
Combining formulas \eqref{eq:reccurence} and \eqref{eq:ratio}, we obtain an irrational ratio $ k=\sqrt{2}-1 $. This ratio suggests that the Jerusalem square cannot be built from a simple integer grid \cite{BairdCube, BairdSquare}. 

\begin{figure}[H]
\centering
\includegraphics[width=0.25\linewidth]{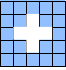}
\caption{An integer grid approximation to the first iteration of the Jerusalem square.}
\label{Fig:L1approx}
\end{figure}

However, a naive method is to consider a $ 5\times 5  $ square, and then remove a cross which consists of five unit squares as shown in \Figs{Fig:L1approx}.

We shall see in the next section that it is in fact a nice approximation of the actual fractal for the corresponding iteration.  Indeed, the present paper is motivated by formulas \eqref{eq:reccurence}, \eqref{eq:ratio} and  the observation of \Figs{Fig:L1approx}. On the one hand, notice  the similarities between the recurrence relation \eqref{eq:reccurence} and the definition of the Pell numbers \seqnum{A000129}, and on the other hand, notice that the side lengths of the squares (of rank $ +1 $ and $ +2 $, as well as the original) in \Figs{Fig:L1approx} are exactly $ 1 $, $ 2 $ and $ 5 $, some of the first few terms of the Pell numbers. 

\section{Pell numbers in action}
Firstly, recall that the Pell numbers are defined as the sequence of integers
\begin{equation}
p_0  = 0,\ p_1 = 1,\ \textnormal{and}\ p_{n} = 2p_{n-1} + p_{n-2}\ \textnormal{for all}\ n \geq 2.
\end{equation}
The first few Pell numbers are  
\begin{equation*}
0, 1, 2, 5, 12, 29, 70, 169, 408, 985, 2378,\ldots\ \mbox{(sequence \seqnum{A000129} in the OEIS \cite{Sloane})}.
\end{equation*}
It is well-known that the integer ratio $\dfrac{p_n}{p_{n-1}} $ rapidly approach $ 1+\sqrt{2} $ \cite[p.\ 138]{Koshy}.

Now, let us introduce the following informal notation. 

Let $ P_n $ denote a $ p_n\times p_n $ square whose edges result from the alignment of the squares $ P_{n-1} $,  $ P_{n-2} $ and $ P_{n-1} $  as illustrated in the following formula:

$ P_0 :=\varnothing$ (the square of side length $ 0 $), $ P_1 :=\protect\includegraphics[width=.03\linewidth,valign=c]{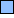}$ (the unit square) and
\begin{equation}\label{Eq:DefinitionPn}
P_n:=\left[\begin{matrix}
P_{n-1}&P_{n-2}&P_{n-1}\\
P_{n-2}& &P_{n-2}\\
P_{n-1}&P_{n-2}&P_{n-1}\\
\end{matrix}\right]=\protect\includegraphics[width=.31\linewidth,valign=c]{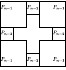}\ \textnormal{for all}\ n \geq 2.
\end{equation}
The blank entry in the matrix representation in \eqref{Eq:DefinitionPn} is there to indicate the cross removal. For example, for $ n=2,3,4,5,6 $ we have
\begin{equation*}
P_2 =\left[\begin{matrix}
P_{1}&P_{0}&P_{1}\\
P_{0}& &P_{0}\\
P_{1}&P_{0}&P_{1}\\
\end{matrix}\right]=\protect\includegraphics[width=.05\linewidth,valign=c]{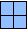},\ 
P_3 =\left[\begin{matrix}
P_{2}&P_{1}&P_{2}\\
P_{1}& &P_{1}\\
P_{2}&P_{1}&P_{2}\\
\end{matrix}\right]=\protect\includegraphics[width=.11\linewidth,valign=c]{level1},
\end{equation*} 
\begin{equation*}
P_4 =\left[\begin{matrix}
P_{3}&P_{2}&P_{3}\\
P_{2}& &P_{2}\\
P_{3}&P_{2}&P_{3}\\
\end{matrix}\right]=\protect\includegraphics[width=.15\linewidth,valign=c]{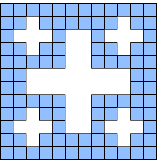},\  
P_5 =\left[\begin{matrix}
P_{4}&P_{3}&P_{4}\\
P_{3}& &P_{3}\\
P_{4}&P_{3}&P_{4}\\
\end{matrix}\right]=\protect\includegraphics[width=.3\linewidth,valign=c]{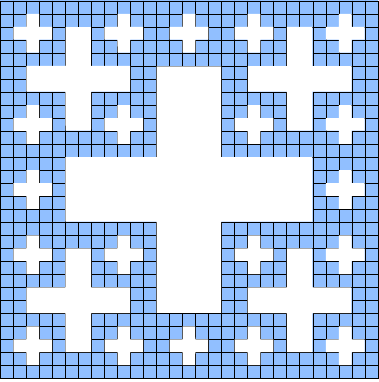},
\end{equation*}
\begin{equation*}
P_6 =\left[\begin{matrix}
P_{5}&P_{4}&P_{5}\\
P_{4}& &P_{4}\\
P_{5}&P_{4}&P_{5}\\
\end{matrix}\right]=\protect\includegraphics[width=.58\linewidth,valign=c]{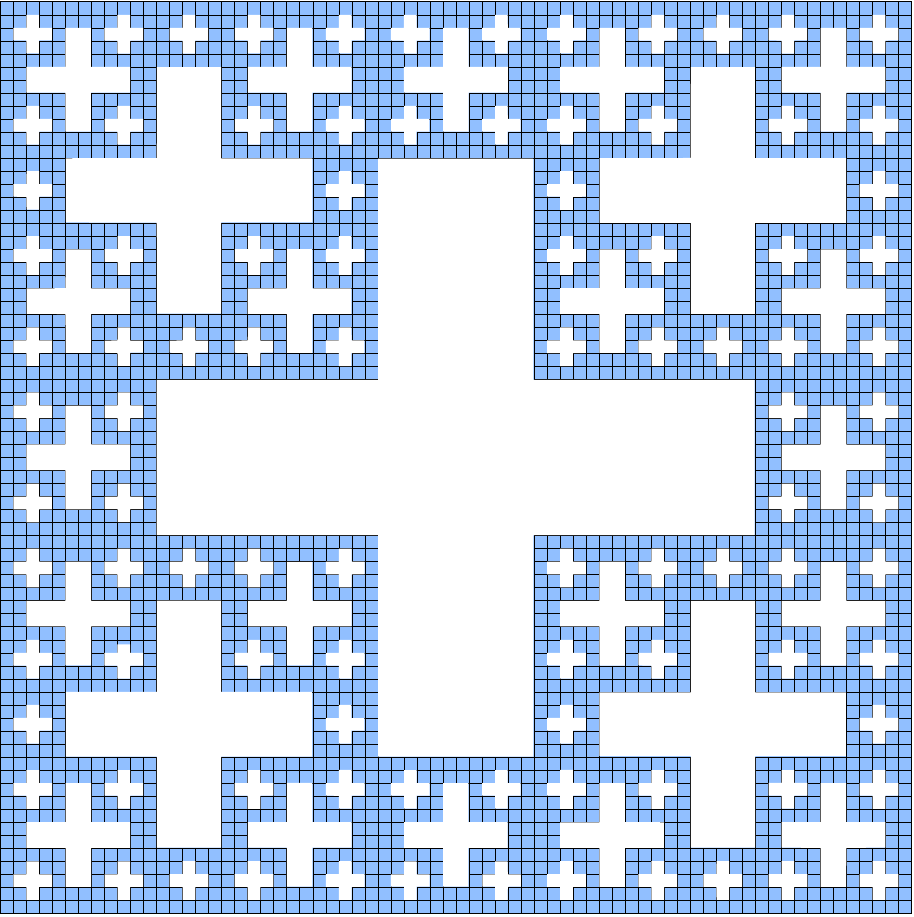}.
\end{equation*}

Since the square $ P_n $ is of length $ p_n $, then by the property of the Pell numbers, the ratios $\dfrac{p_{n-1}}{p_{n}}$ and $\dfrac{p_{n-2}}{p_{n-1}}$ give a good approximation to the Jerusalem square ratio $\sqrt{2}-1 $ when $ n $ is sufficiently large. 

\begin{figure}[H]
\centering
\hspace*{\fill}
\subfigure[$ n=2 $]{\includegraphics[width=0.1\linewidth]{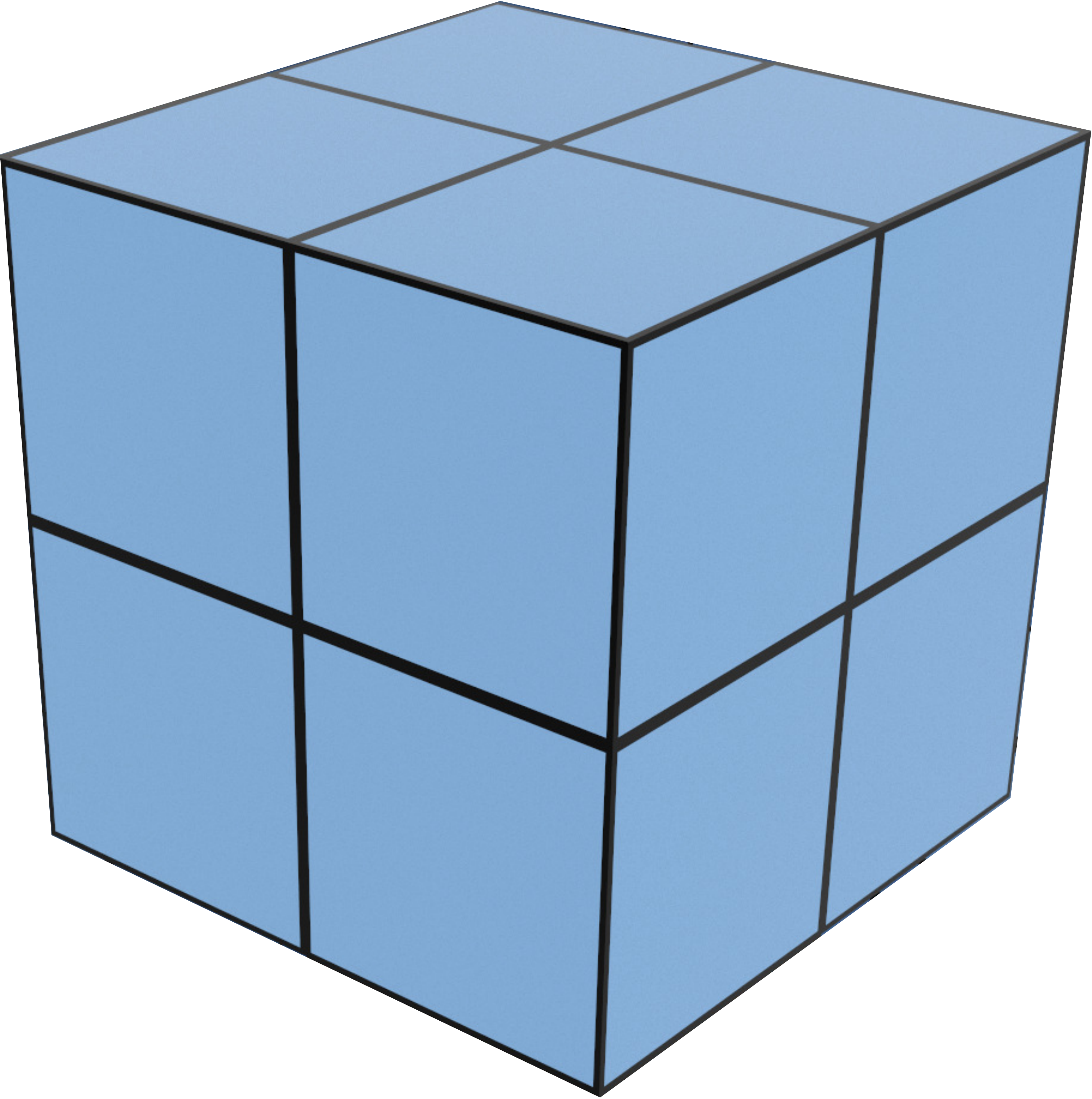}}\hfill%
\subfigure[$ n=3 $]{\includegraphics[width=0.3\linewidth]{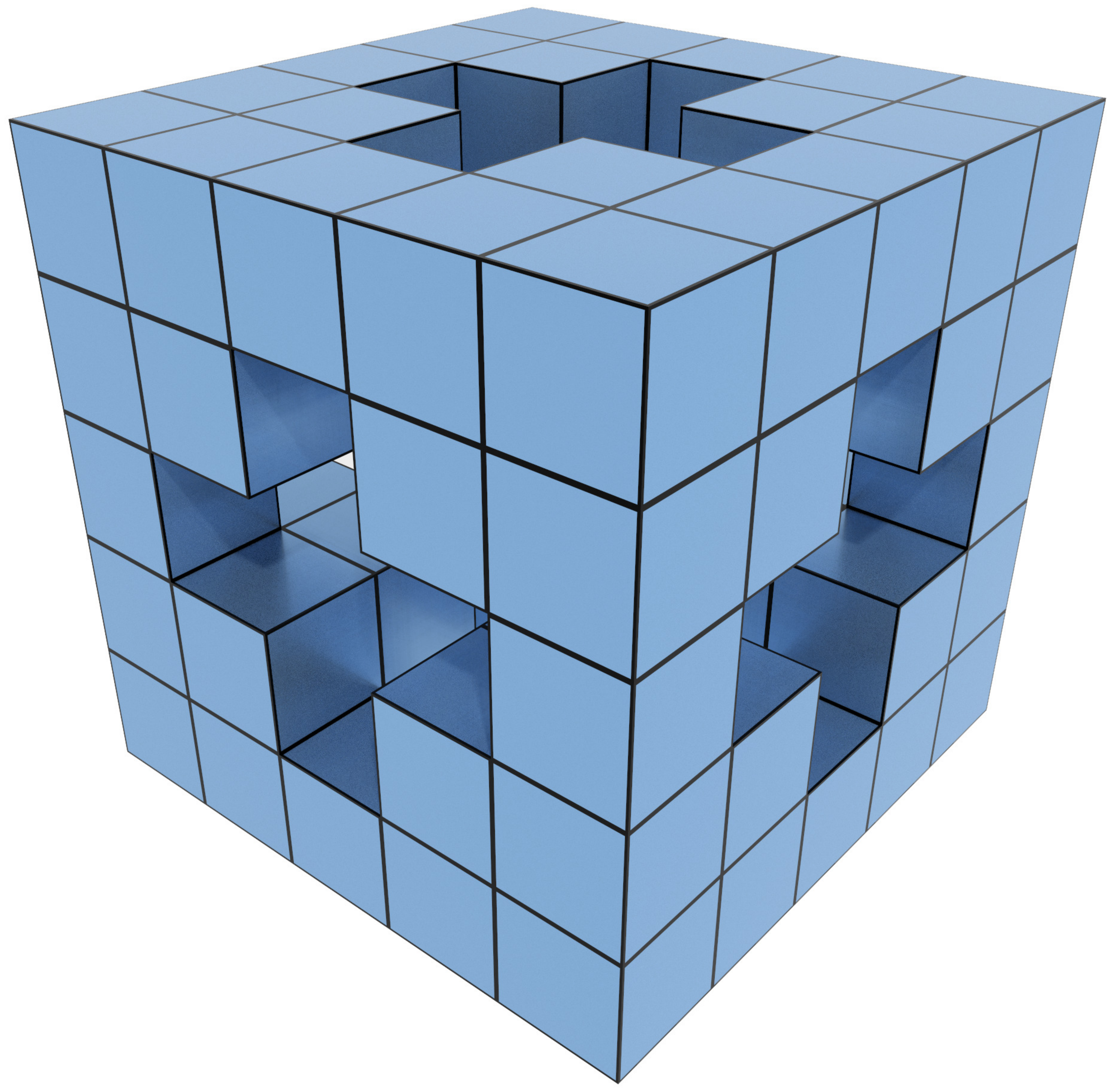}}\hfill%
\subfigure[$ n=4 $]{\includegraphics[width=0.45\linewidth]{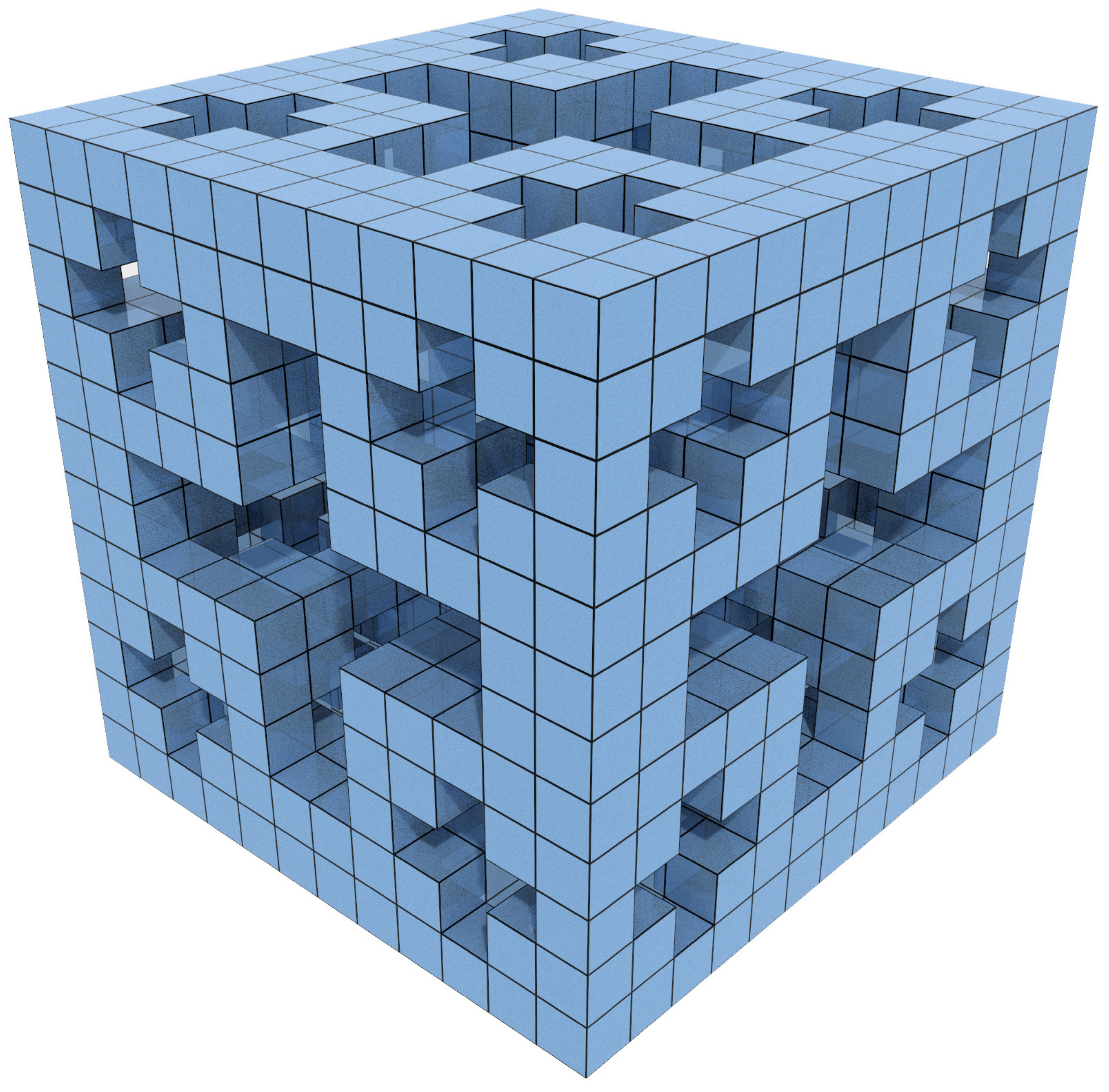}}\hfill%
\subfigure[$ n=5$]{\includegraphics[width=0.4\linewidth]{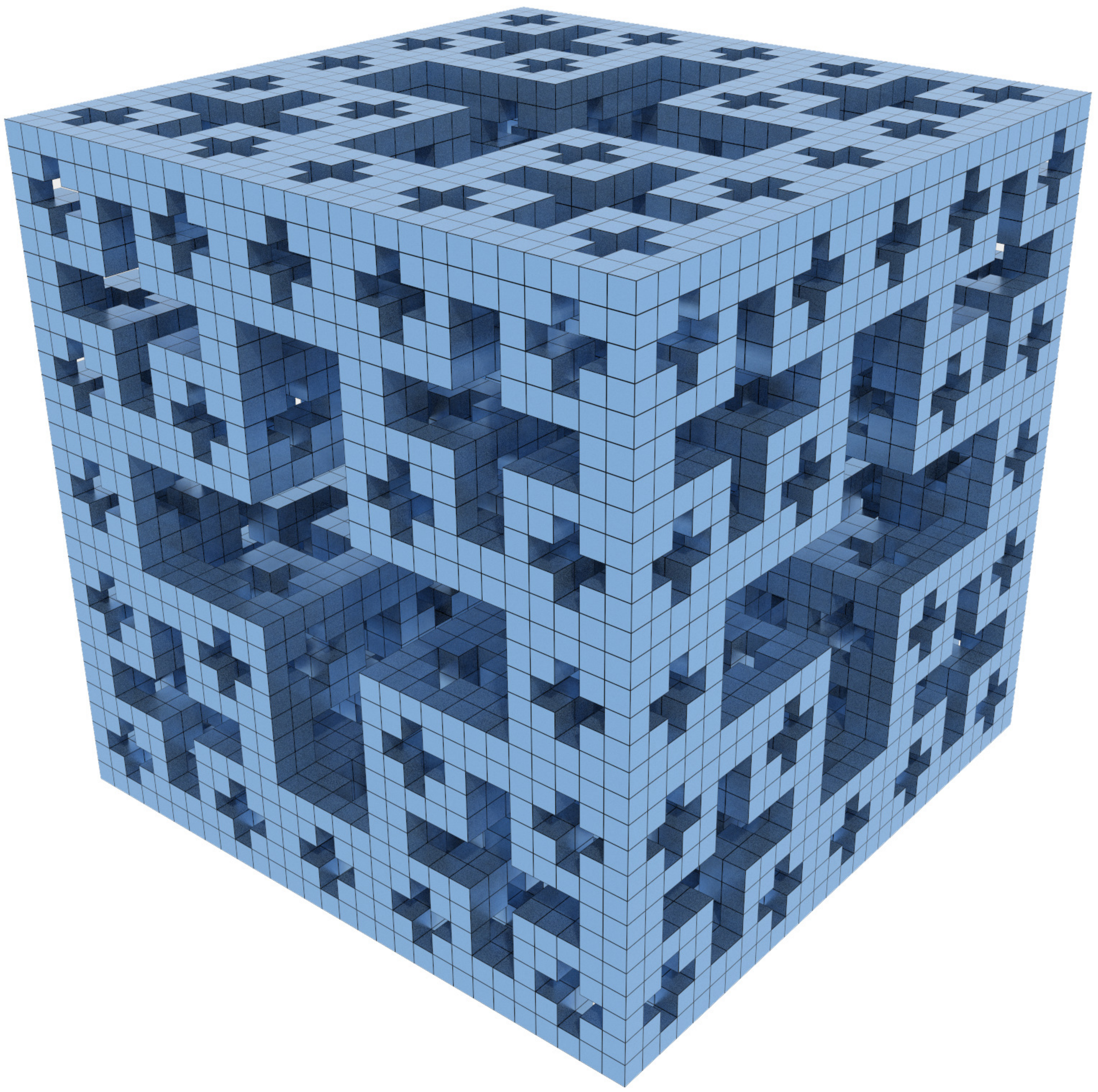}}\hfill%
\subfigure[$ n=6$]{\includegraphics[width=0.575\linewidth]{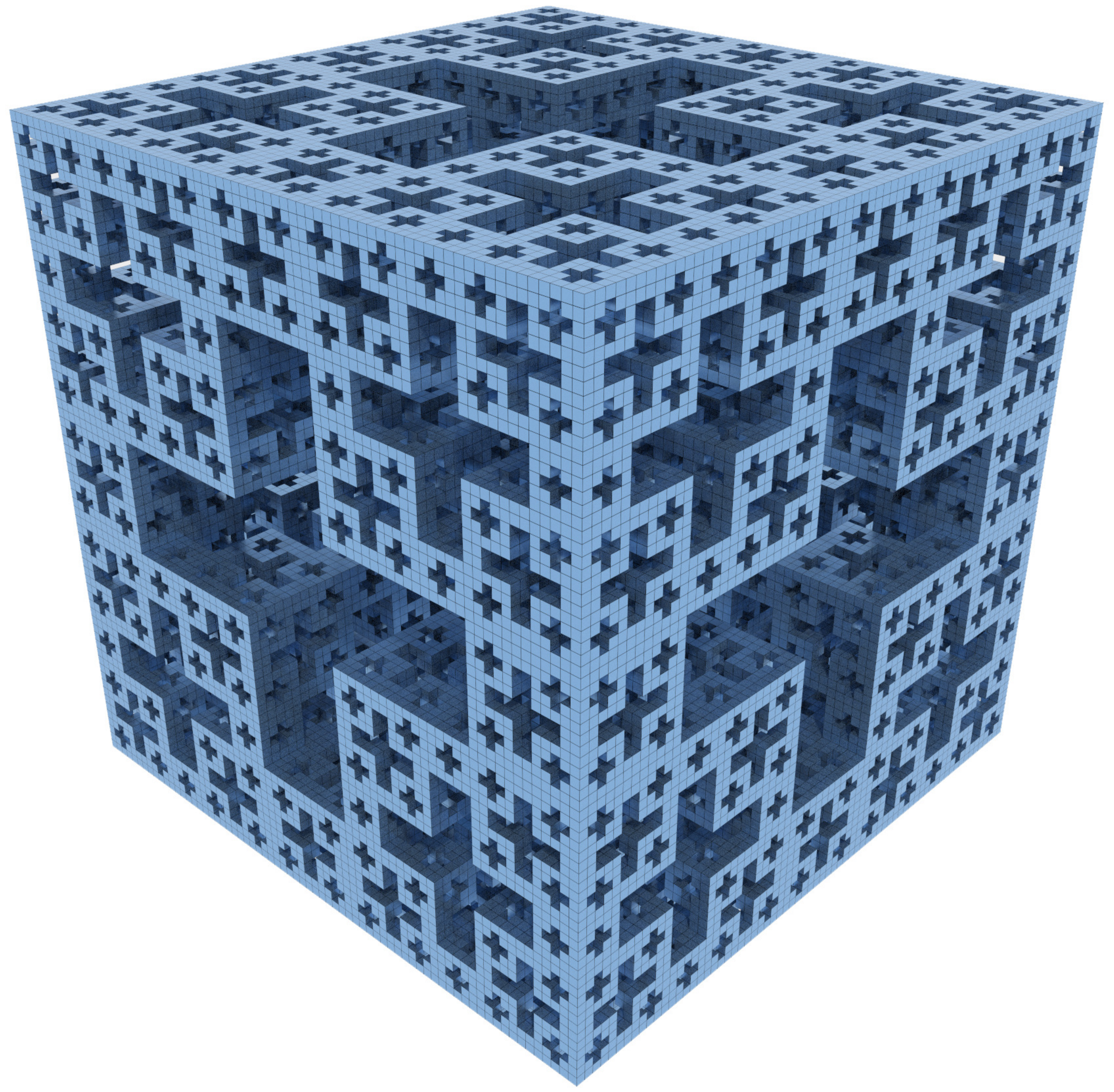}}
\hspace*{\fill}
\caption{Some iterations of the  approximation to the Jerusalem cube.}
\label{Fig:JerusalemCube}
\end{figure}

We can extend this method to the Jerusalem cube \cite{BairdCube}, and consider a construction of this three-dimensional case with the popular business card cube \cite[p.\ 152]{Hull}. For example, we see in \Figs{Fig:JerusalemCube} the corresponding iterations for $ n=2,3,4,5,6 $.

\end{document}